# Asymptotic of some sums

VICTOR VOLFSON

ABSTRACT.

The paper compares the asymptotic of the expressions $\frac{1}{x}\sum_{n\leq x} f(n)$ and $\sum_{n\leq x}\frac{f(n)}{n}$, $\frac{1}{x}\sum_{n\leq x} f(p)$ and $\sum_{n\leq x}\frac{f(p)}{p}$. The asymptotic of sums $\sum_{p\leq x}\frac{f(n)}{n}$ and $\sum_{p\leq x}\frac{f(p)}{p}$ ($n, p$ - respectively, positive and prime numbers) are determined if the asymptotic of sums are known, respectively: $\sum_{n\leq x} f(n), \sum_{p\leq x} f(p)$.

1. INTRODUCTION

The goal of this paper is to determine the asymptotic of sums $\sum_{p\leq x}\frac{f(n)}{n}$ and $\sum_{p\leq x}\frac{f(p)}{p}$ ($n, p$ - respectively, a natural and prime number), if the asymptotic of sums are known, respectively: $\sum_{n\leq x} f(n), \sum_{p\leq x} f(p)$.

We note that in this way it is also possible to determine the asymptotic of summing arithmetic functions $\sum_{n\leq x}\frac{f(n)}{n^k}$ and $\sum_{n\leq x}\frac{f(p)}{p^k}$, where $k$ is a positive integer.

The principal limitation in determining these asymptotic is that the corresponding integral should be an elementary function. This will become clear from the further presentation.

---

## 2. ABEL SUMMATION FORMULA

Let $a_n$ is a sequence of real or complex numbers and $f(x)$ is a function continuously differentiable on the ray $[1, x)$. Then, if we base on [1]:

$$\sum_{1 \le n \le x} a_n f(n) = A(x) f(x) - \int_1^x A(u) f'(u) du, \qquad (2.1)$$

where $A(x) = \sum_{1 \le n \le x} f(n)$.

We obtain the formula in the particular case if $f(x) = 1/x$ in (2.1):

$$\sum_{n \le x} \frac{f(n)}{n} = \frac{1}{x} \sum_{n \le x} f(n) + \int_1^x \sum_{n \le x} f(n) \frac{dx}{x^2}. \qquad (2.2)$$

Suppose we know the meaning $B(x) = \frac{1}{x} \sum_{n \le x} f(n)$.

Then the formula (2.2) can be written in the form:

$$\sum_{n \le x} \frac{f(n)}{n} = B(x) + \int_1^x \frac{B(t) dt}{t}. \qquad (2.3)$$

Suppose that an asymptotic estimate is known $\sum_{n \le x} f(n)$, it is easily determined $B(x) = \frac{1}{x} \sum_{n \le x} f(n)$ and the integral $\int_1^x \frac{B(t) dt}{t}$ if it exists in elementary functions.

Suppose that an asymptotic estimate is known $\sum_{p \le x} f(p)$, then it is easily determined $B(x) = \frac{1}{x} \sum_{n \le x} f(p)$ and also it is easily determined the integral $\int_1^x \frac{B(t) dt}{t}$ if it exists in elementary functions.

Therefore, the formula is true similar to (2.3):

$$\sum_{p \le x} \frac{f(p)}{p} = B(x) + \int_1^x \frac{B(t) dt}{t}, \qquad (2.4)$$

where $p$ is a prime number.

Based on (2.3), it is clear that if the integral on the right has an asymptotic upper bound that is superior $B(x)$, then the upper asymptotic bound $\sum_{n \leq x} \frac{f(n)}{n}$ exceeds the upper one $\frac{1}{x}\sum_{n \leq x} f(n)$, if not, then the upper asymptotic estimates $\frac{1}{x}\sum_{n \leq x} f(n)$ and $\sum_{n \leq x} \frac{f(n)}{n}$ coincides.

Having in mind (2.4), it is true that if the integral on the right has an asymptotic upper bound that is superior $B(x)$, then the upper asymptotic estimate $\sum_{n \leq x} \frac{f(p)}{p}$ exceeds the upper bound $\frac{1}{x}\sum_{n \leq x} f(p)$, if not, then the indicated asymptotic estimates coincide.

3. COMPARING ASYMPTOTIC EXPRESSIONS $\frac{1}{x}\sum_{n \leq x} f(n)$ AND $\sum_{n \leq x} \frac{f(n)}{n}$, $\frac{1}{x}\sum_{n \leq x} f(p)$ AND $\sum_{n \leq x} \frac{f(p)}{p}$

Let's start the comparison from the simplest function $f(n) = 1$.

Naturally meaning:

$$\frac{1}{x}\sum_{n \leq x} 1 = 1. \tag{3.1}$$

Based on (2.3) and (3.1) we get:

$$\sum_{n \leq x} \frac{1}{n} = 1 + \int_1^x \frac{dt}{t} = 1 + \log(x) = O(\log(x)). \tag{3.2}$$

The asymptotic upper bound (3.2) has a larger order than (3.1), since the integral has an upper asymptotic upper bound of order than (3.1).

Based on the Euler-Maclaurin formula [2] for the function $f(n) = \log^k(n) + O(\log^{k-1} n)$ ($k > 0$), we obtain the following estimate:

$$\frac{1}{x}\sum_{n \leq x} \log^k n = \log^k x + O(\log^{k-1} x). \tag{3.3}$$

Having in mind the Euler-Maclaurin formula, we obtain the following estimate:

$$\sum_{n \leq x} \frac{\log^k n}{n} = \int_1^x \frac{\log^k t}{t} dt + c + O(\frac{\log^k x}{x}) = \frac{\log^{k+1} x}{k+1} + c + O(\frac{\log^k x}{x}). \qquad (3.4)$$

The asymptotic estimate (3.4) has a larger order than (3.3), since the integral in (2.3) has an asymptotic estimate of a higher order than (3.3).

Based on the Euler-Maclaurin formula for the function $f(n) = n^m + O(n^{m-1})$ ($m > 0$), we obtain the following estimate:

$$\frac{1}{x}\sum_{n \leq x} n^m = \frac{1}{x}[\int_1^x t^m dt + O(x^m)] = \frac{x^m}{m+1} + O(x^{m-1}) = O(x^m). \qquad (3.5)$$

Having in mind the Euler-Maclaurin formula, we obtain the following estimate:

$$\sum_{n \leq x} \frac{n^m}{n} = \frac{x^m}{m+1} + O(x^{m-1}) + \frac{1}{m+1}\int_1^x t^{m-1} dt = (1+\frac{1}{m})\frac{x^m}{m+1} + O(x^{m-1}) = O(x^m). \qquad (3.6)$$

The asymptotic upper bounds (3.5) and (3.6) coincide, since the integral in (2.3) has the same asymptotic upper bound as (3.5).

Now we consider Chebyshev functions - $\Psi(x), \theta(x)$. It is known [3] that they have asymptotic:

$$\Psi(x) = x + o(x), \theta(x) = x + o(x). \qquad (3.7)$$

Based on (3.7), we obtain the average values of the terms of these functions:

$$\frac{\Psi(x)}{x} = 1 + o(1), \frac{\theta(x)}{x} = 1 + o(1). \qquad (3.8)$$

Let us consider, for example, $\theta(x) = \sum_{p \leq x} \log(p) = x + o(x)$, where $p$ is a prime number.

Thus, the following asymptotic estimate holds:

$$\frac{1}{x}\sum_{p \leq x} \log(p) = \frac{\theta(x)}{x} = 1 + o(1). \qquad (3.9)$$

Having in mind (2.4) and (3.9) we get:

$$\sum_{p \leq x} \frac{\log(p)}{p} = 1 + o(1) + \int_1^x \frac{dt}{t} + O(\log(x)) = O(\log(x)). \tag{3.10}$$

The asymptotic upper bound (3.10) is of order greater than (3.9), since the integral has an asymptotic estimate of a higher order than (3.9).

The following upper bound for Möbius function is known [4]:

$$|\sum_{n \leq x} \mu(n)| \leq \frac{c_2 x}{\log^2 x}, \tag{3.11}$$

where $\mu(n)$ is Möbius function and $c_2 = 362,7$.

Based on (3.11), the following asymptotic estimate from above is satisfied:

$$\frac{1}{x} |\sum_{n \leq x} \mu(n)| = O(\frac{1}{\log^2 x}). \tag{3.12}$$

Then, having in mind (2.3) and (3.12) we get:

$$|\sum_{n \leq x} \frac{\mu(n)}{n}| \leq \frac{c_1}{\log^2 x} + \int_1^x \frac{c_1 dt}{t \log^2 t} = \frac{c_1}{\log^2 x} + \frac{c_2}{\log^3 x},$$

so:

$$|\sum_{n \leq x} \frac{\mu(n)}{n}| = O(\frac{1}{\log^2 x}). \tag{3.13}$$

The asymptotic upper bounds (3.12) and (3.13) coincide, since the integral has an asymptotic upper bound of a smaller order than (3.12).

Considering that:

$$\sum_{n \leq x} \frac{\mu(n)}{n} \leq |\sum_{n \leq x} \frac{\mu(n)}{n}|.$$

we obtain the following asymptotic estimate from above:

$$\sum_{n \leq x} \frac{\mu(n)}{n} = O(\frac{1}{\log^2 x}) = o(1). \tag{3.14}$$

# 4. ASYMPTOTICS OF SUMMATORY ARITHMETIC FUNCTIONS OF THE FORM $\sum_{n \leq x} \frac{f(n)}{n}$ AND $\sum_{n \leq x} \frac{f(p)}{p}$

Statement 1

Let the asymptotic estimate is $B(x) = \frac{1}{x}\sum_{n \leq x} f(n) = \frac{1}{\log(x)} + O(\frac{1}{\log^2 x})$, then the following asymptotic estimate is:

$$\sum_{2 \leq n \leq x} \frac{f(n)}{n} = \log\log(x) + O(1). \tag{4.1}$$

Proof

Based on (2.3), we obtain in this case:

$$\sum_{2 \leq n \leq x} \frac{f(n)}{n} = \frac{1}{\log(x)} + O(\frac{1}{\log^2 x}) + \int_2^x \frac{dt}{t\log(t)} + O(1) = \log\log(x) + O(1).$$

Therefore, an asymptotic upper bound (4.1) holds.

As an example, let us consider an arithmetic function of the number of primes that do not exceed the value $x$ - $\pi(x) = \sum_{p \leq x} 1$. Based on the asymptotic law of primes, the following asymptotic is performed:

$$\frac{\pi(x)}{x} = \frac{\sum_{p \leq x} 1}{x} = \frac{1}{\log(x)} + O(\frac{1}{\log^2 x}), \tag{4.2}$$

which corresponds to statement 1.

Therefore, having in mind (4.1), (4.2), the following asymptotic is true:

$$\sum_{p \leq x} \frac{1}{p} = \log\log(x) + O(1). \tag{4.3}$$

Statement 2

Let asymptotic equality hold: $B(x) = c_2 x^m + O(x^{m-1})$, then the following asymptotic estimate from above is true:

$$\sum_{n\leq x}\frac{f(n)}{n}=O(x^m). \qquad (4.4)$$

Based on (3.6) the proof is true.

As an example, we consider the function of the positive divisors of a natural number $n$ - $\sigma(n)$. It is known [5] that:

$$\sum_{n\leq x}\sigma^2(n)=c_3 x^3+O(x^{7/3}\log^3 x).$$

Therefore, the condition of statement 2 is satisfied:

$$\frac{1}{x}\sum_{n\leq x}\sigma^2(n)=c_3 x^2+O(x^{4/3}\log^3 x)=c_3 x^2+O(x). \qquad (4.5)$$

Having in mind (4.4) and (4.5), we obtain the following asymptotic upper bound:

$$\sum_{n\leq x}\frac{\sigma^2(n)}{n}=O(x^2). \qquad (4.6)$$

Based on (2.3) and (4.5), we can refine the estimate (4.6):

$$\sum_{n\leq x}\frac{\sigma^2(n)}{n}=c_4 x^2+O(x^{4/3}\log^3 x). \qquad (4.7)$$

Statement 3

Let us the asymptotic equality hold: $B(x)=\log^k(x)+O(\log^{k-1} x)$, then the following asymptotic estimate from above is true:

$$\sum_{n\leq x}\frac{f(n)}{n}=O(\log^{k+1} x). \qquad (4.8)$$

The proof follows from formula (3.4).

As an example, we consider the well-known asymptotic of the following function [6]:

$$\sum_{n\leq x}\tau(n)=x\log(x)+c_6 x+o(x), \qquad (4.9)$$

where the function $\tau(n)$ is the number of positive divisors of a natural number $n$.

Having in mind (4.9) we find the asymptotic:

$$B(x) = \frac{1}{x} \sum_{n \leq x} \tau(n) = \log(x) + O(1), \qquad (4.10)$$

which is subject to Statement 3.

Therefore, based on (4.8) and (4.9), we obtain the following asymptotic estimate from above:

$$\sum_{n \leq x} \frac{\tau(n)}{n} = O(\log^2 x). \qquad (4.11)$$

Having in mind (2.3) and (4.10), we can refine the estimate (4.11):

$$\sum_{n \leq x} \frac{\tau(n)}{n} = \log(x) + O(1) + \int_1^x \frac{\log(t)dt}{t} + O(\log(x)) = \frac{\log^2 x}{2} + O(\log(x)). \qquad (4.12)$$

## 5. CONCLUSION AND SUGGESTIONS FOR FURTHER WORK

The next article will continue to study the behavior of some sums.

## 6. ACKNOWLEDGEMENTS

Thanks to everyone who has contributed to the discussion of this paper. I am grateful to everyone who expressed their suggestions and comments in the course of this work.